\input amstex
\documentstyle{amsppt}
%
\catcode`@=11
\redefine\output@{%
  \def\break{\penalty-\@M}\let\par\endgraf
  \ifodd\pageno\global\hoffset=105pt\else\global\hoffset=8pt\fi  
  \shipout\vbox{%
    \ifplain@
      \let\makeheadline\relax \let\makefootline\relax
    \else
      \iffirstpage@ \global\firstpage@false
        \let\rightheadline\frheadline
        \let\leftheadline\flheadline
      \else
        \ifrunheads@ 
        \else \let\makeheadline\relax
        \fi
      \fi
    \fi
    \makeheadline \pagebody \makefootline}%
  \advancepageno \ifnum\outputpenalty>-\@MM\else\dosupereject\fi
}
\def\Beta{\mathchar"0\hexnumber@\rmfam 42}
\catcode`\@=\active
\nopagenumbers
\chardef\textvolna='176

\chardef\bigalpha='013

\def\const{\operatorname{const}}

\def\compos{\,\raise 1pt\hbox{$\sssize\circ$} \,}
\def\Span{\operatorname{Span}}
\def\msum#1{\operatornamewithlimits{\sum^#1\!{\ssize\ldots}\!\sum^#1}}

\def\blue#1{#1}

\catcode`#=11\def\diez{#}\catcode`#=6
\catcode`&=11\catcode`&=4
\catcode`_=11\def\podcherkivanie{_}\catcode`_=8
\catcode`~=11\catcode`~=\active
\def\mycite#1{\cite{\blue{#1}}\immediate\special{ps:
     ShrHPSdict begin /ShrBORDERthickness 0 def}}

\def\mytag#1{%
    \tag#1}
\def\mythetag#1{\thetag{\blue{#1}}\immediate\special{ps:
     ShrHPSdict begin /ShrBORDERthickness 0 def}}
\def\myrefno#1{\no#1}
\def\myhref#1#2{\blue{#2}\immediate\special{ps:
     ShrHPSdict begin /ShrBORDERthickness 0 def}}

\def\mytheorem#1{\csname proclaim\endcsname{Theorem #1}}
\def\mytheoremwithtitle#1#2{\csname proclaim\endcsname{Theorem #1#2}}
\def\mythetheorem#1{\blue{#1}\immediate\special{ps:
     ShrHPSdict begin /ShrBORDERthickness 0 def}}
\def\mylemma#1{\csname proclaim\endcsname{Lemma #1}}
\def\mylemmawithtitle#1#2{\csname proclaim\endcsname{Lemma #1#2}}
\def\mythelemma#1{\blue{#1}\immediate\special{ps:
     ShrHPSdict begin /ShrBORDERthickness 0 def}}
\def\mycorollary#1{\csname proclaim\endcsname{Corollary #1}}

\def\mydefinition#1{\definition{Definition #1}}
\def\mythedefinition#1{\blue{#1}\immediate\special{ps:
     ShrHPSdict begin /ShrBORDERthickness 0 def}}
\def\myconjecture#1{\csname proclaim\endcsname{Conjecture #1}}
\def\myconjecturewithtitle#1#2{\csname proclaim\endcsname{Conjecture #1#2}}

\def\myproblem#1{\csname proclaim\endcsname{Problem #1}}
\def\myproblemwithtitle#1#2{\csname proclaim\endcsname{Problem #1#2}}


\pagewidth{360pt}
\pageheight{606pt}
\topmatter
\title
Clusters of exponential functions in the space of square integrable functions.
\endtitle
\rightheadtext{Clusters of exponents \dots}
\author
Ruslan Sharipov
\endauthor
\address Bashkir State University, 32 Zaki Validi street, 450074 Ufa, Russia
\endaddress
\email\myhref{mailto:r-sharipov\@mail.ru}{r-sharipov\@mail.ru}
\endemail
\abstract
     Finite dimensional subspaces spanned by exponential functions in the space
of square integrable functions on a finite interval of the real line are considered.
Their limiting positions are studied and described in terms of expo-polynomials. 
\endabstract
\subjclassyear{2000}
\subjclass 42C15, 46C07 \endsubjclass
\endtopmatter
\TagsOnRight
\document

\head
1. Introduction.
\endhead
     Exponential functions of the form $e^{\kern 0.5pt\lambda\,x}$ naturally 
arise in Fourier series. More general series of exponential functions were 
studied by A. F. Leontiev (see \mycite{1}) and his school. In this paper we 
consider finite sequences of exponential functions
$$
\hskip -2em
e^{\kern 0.5pt\lambda_1x},\,\ldots,\,e^{\kern 0.5pt\lambda_nx},
\mytag{1.1}
$$
where $\lambda_1,\,\ldots,\,\lambda_n$ are distinct complex numbers, i\.\,e\.
$\lambda_{\kern 1pt i}\neq\lambda_j$. These numbers are called the spectrum of the
sequence \mythetag{1.1}.\par
     The exponential functions \mythetag{1.1} are treated as elements of the 
space of square integrable functions $L^2([a,b])$, where $-\infty<a<b<+\infty$.
Without loss of generality we can take $a=-\pi$ and $b=+\pi$ like in Fourier
analysis. The functions \mythetag{1.1} span an $n$-dimensional subspace in
$L^2([-\pi,+\pi])$:
$$
\hskip -2em
L=\Span(e^{\kern 0.5pt\lambda_1x},\,\ldots,e^{\kern 0.5pt\lambda_nx}).
\mytag{1.2}
$$ 
The goal of this paper is to describe the behavior of the subspace \mythetag{1.2}
when the numbers $\lambda_1,\,\ldots,\,\lambda_n$ subdivide into $m$ clusters and 
tend to some distinct limit values $\Lambda_1,\,\ldots,\,\Lambda_n$ common within 
each cluster. Therefore below we use the following double index notation for the 
numbers $\lambda_1,\,\ldots,\,\lambda_n$:
$$
\hskip -2em
\lambda_{\kern 0.5pt ij}\text{, \ where \ }
i=1,\,\ldots,m\text{\ \ and \ }j=1,\,\ldots,k_i. 
\mytag{1.3}
$$
The numbers $k_1,\ldots,\,k_m$ in \mythetag{1.3} are called multiplicities of
clusters. Their sum is equal to the total number of lambdas in \mythetag{1.2}:
$$
\hskip -2em
k_1+\ldots+k_m=n.
\mytag{1.4}
$$
Due to \mythetag{1.3} and \mythetag{1.4} we write \mythetag{1.2} as 
$$
L=\Span(\{e^{\kern 0.5pt\lambda_{ij}x}\}^{i=1,\,\ldots,\,m}_{j=1,\,\ldots,\,k_i})
$$ 
or as $L=\Span(\{e^{\kern 0.5pt\lambda_{ij}x}\})$ for short. As it was said above,
we assume that
$$
\hskip -2em
\lambda_{ij}\to\Lambda_i
\mytag{1.5}
$$
in a sequence of samplings or in a continuous process. It is convenient to denote
$$
\hskip -2em
\varepsilon=\max\{|\lambda_{ij}-\Lambda_i|\}^{i=1,\,\ldots,\,m}_{j=1,\,\ldots,\,k_i}.
\mytag{1.6}
$$
Then we can write the formula \mythetag{1.5} in the following way:
$$
\hskip -2em
\lambda_{ij}\to\Lambda_i\text{ \ as \ }\varepsilon\to 0.
\mytag{1.7}
$$
\head
2. Convergence of subspaces in a Hilbert space.
\endhead
     Let $\Cal H$ be a Hilbert space (see \mycite{2}). The space of square
integrable functions $L^2([-\pi,+\pi])$ is an example of a Hilbert space. 
\mydefinition{2.1} A sequence $L_q$ of $n$-dimensional subspaces of a Hilbert 
space $\Cal H$ is said to converge to an $n$-dimensional subspace $M$ if there
are some bases $\bold e_{1q},\,\ldots,\,\bold e_{nq}$ in $L_q$ and there is
some basis $\bold e_{1},\,\ldots,\,\bold e_{n}$ in $M$ such that 
$$
\bold e_{iq}\to \bold e_i\text{\ \ as \ }q\to\infty
$$
in the sense of the norm of the Hilbert space $\Cal H$. 
\enddefinition
     The definition~\mythedefinition{2.1} can be reformulated in order to apply 
to the case of continuous parametric sets of subspaces. 
\mydefinition{2.2} Let $L_\varepsilon$ be a parametric set of $n$-dimensional 
subspaces of a Hilbert space $\Cal H$.  It is said to converge to an $n$-dimensional 
subspace $M$ as $\varepsilon\to 0$ if there are some bases $\bold e_{1\varepsilon},
\,\ldots,\,\bold e_{n\varepsilon}$ in $L_\varepsilon$ and there is a basis 
$\bold e_{1},\,\ldots,\,\bold e_{n}$ in $M$ such that 
$$
\bold e_{i\varepsilon}\to \bold e_i\text{\ \ as \ }\varepsilon\to 0
$$
in the sense of the norm of the Hilbert space $\Cal H$. 
\enddefinition
\head
3. Taylor expansions of exponential functions. 
\endhead
     This section is a preliminary one. Assume for a while that we have only one
cluster (i\.\,e\. $m=1$) with $\Lambda_1=0$. Then we can use the initial notations
$\lambda_1,\,\ldots,\,\lambda_n$ for lambdas and write the formula \mythetag{1.5} as
$\lambda_{\kern 1pt i}\to 0$. The exponential functions \mythetag{1.2} have the following 
Taylor expansions:
$$
\hskip -2em
\aligned
&e^{\kern 0.5pt\lambda_1x}=1+\lambda_1\,x+\ldots
+\frac{\lambda_1^{n-1}\,x^{n-1}}{(n-1)!}+\ldots
=\sum^\infty_{q=0}\frac{\lambda_1^q\,x^q}{q!},\\
\vspace{-1ex}
&.\ .\ .\ .\ .\ .\ .\ .\ .\ .\ .\ .\ .\ .\ .\ .\ .\ .\ .\ .\ .\ .\ .\ .\ 
.\ .\ .\ .\ .\ .\ .\ .\ .\ .\ .\ .\ .\ .\ 
\\
\vspace{-0.4ex}
&e^{\kern 0.5pt\lambda_nx}=1+\lambda_n\,x+\ldots
+\frac{\lambda_n^{n-1}x^{n-1}}{(n-1)!}+\ldots
=\sum^\infty_{q=0}\frac{\lambda_n^qx^q}{q!}. 
\endaligned
\mytag{3.1}
$$
Initial parts of the power series \mythetag{3.1} define the polynomials
$$
\hskip -2em
p_i(x)=1+\lambda_{\kern 1pt i}\,x+\ldots+\frac{\lambda_{\kern 1pt i}^{n-1}\,x^{n-1}}{({n-1})!}
\text{, \ where \ }i=1,\,\ldots,\,n.
\mytag{3.2}
$$
Using the polynomials \mythetag{3.2}, we define the following equation for 
$a_1,\,\ldots,\,a_n$:
$$
\hskip -2em
a_1\,p_1(x)+\ldots+a_n\,p_n(x)=\frac{x^{n-1}}{(n-1)!}. 
\mytag{3.3}
$$
The polynomial equation \mythetag{3.3} is equivalent to a matrix equation for 
$a_1,\,\ldots,\,a_n$:
$$
\hskip -2em
\Vmatrix
1 & 1 &\hdots & 1\\
\lambda_1 & \lambda_2 &\hdots & \lambda_n\\
\vspace{0.5ex}
\lambda_1^2 & \lambda_2^2 &\hdots & \lambda_n^2\\
\vdots & \vdots & \ddots & \vdots\\ 
\vspace{0.5ex}
\lambda_1^{n-1} & \lambda_2^{n-1} &\hdots & \lambda_n^{n-1}
\endVmatrix
\cdot
\Vmatrix
a_1\\ 
\vspace{0.5ex}
a_2\\
\vspace{0.5ex}
a_3\\ \vdots\\ a_n
\endVmatrix
=\Vmatrix
0\\ 
\vspace{0.5ex}
0\\
\vspace{0.5ex}
0\\ \vdots\\ 1
\endVmatrix.
\mytag{3.4}
$$
The matrix in \mythetag{3.4} is the transpose of the Vandermonde 
matrix (see \mycite{3}):
$$
\hskip -2em
W=\Vmatrix
\ 1 & \lambda_1 & \lambda_1^2 & \hdots & \lambda_1^{n-1}\\
\vspace{0.8ex}
\ 1 & \lambda_2 & \lambda_2^2 & \hdots & \lambda_2^{n-1}\\
\vspace{0.8ex}
\ 1 & \lambda_3 & \lambda_3^2 & \hdots & \lambda_3^{n-1}\\
\ \vdots & \vdots & \vdots & \ddots & \vdots\\
\vspace{0.7ex}
\ 1 & \lambda_n & \lambda_n^2 & \hdots & \lambda_n^{n-1}
\endVmatrix.
\mytag{3.5}
$$\par 
The Vandermonde matrix \mythetag{3.5} is non-degenerate for distinct lambdas,
i\.\,e\. if $\lambda_{\kern 1pt i}\neq\lambda_j$. In this case it has the inverse matrix
$U=W^{-1}$ (see \mycite{4}). In order to write the elements of the inverse matrix
$U=W^{-1}$ explicitly we use the polynomials
$$
\hskip -2em
P_q(\lambda)=\frac{\dsize\prod^n_{s\neq q}(\lambda-\lambda_s)}
{\dsize\prod^n_{s\neq q}(\lambda_q-\lambda_s)}
\text{, \ where \ }q=1,\,\ldots,\,n. 
\mytag{3.6}
$$
It is easy to see that the polynomials \mythetag{3.6} obey the equality
$$
\hskip -2em
P_q(\lambda_{\kern 1pt i})=\cases 1 &\text{for \ }q=i,\\
0 &\text{for \ }q\neq i.\endcases
\mytag{3.7}
$$
If we present the polynomials \mythetag{3.6} as the power expansions 
$$
\hskip -2em
P_q(\lambda)=\sum^n_{r=1}U_{rq}\,\lambda^{r-1}=
U_{1q}+U_{2q}\,\lambda+\ldots+U_{nq}\,\lambda^{n-1},  
\mytag{3.8}
$$
then the equality \mythetag{3.7} can be rewritten as 
$$
\hskip -2em
\sum^n_{r=1}\lambda_{\kern 1pt i}^{r-1}\,U_{rq}=\cases 1 &\text{for \ }q=i,\\
0 &\text{for \ }q\neq i.\endcases
\mytag{3.9}
$$
Looking at the matrix \mythetag{3.5}, we see that the equality
\mythetag{3.9} is equivalent to the matrix equality $W\cdot U=1$, where
$U$ is the matrix whose components coincide with the coefficients
in power expansions \mythetag{3.8}:
$$
\hskip -2em
U=\Vmatrix
U_{11} & U_{12} & U_{13} & \hdots & U_{1n}\\
\vspace{0.8ex}
U_{21} & U_{22} & U_{23} & \hdots & U_{2n}\\
\vspace{0.8ex}
U_{31} & U_{32} & U_{33} & \hdots & U_{3n}\\
\ \vdots & \vdots & \vdots & \ddots & \vdots\\
\vspace{0.7ex}
U_{n1} & U_{n2} & U_{n3} & \hdots & U_{nn}\\
\endVmatrix.
\mytag{3.10}
$$
Explicit expressions for the components of the matrix \mythetag{3.10}
are derived from \mythetag{3.6}:
$$
\hskip -2em
U_{rq}=\frac{1}{(r-1)!}\,\frac{d^{\kern 1pt r-1}P_q(\lambda)}{d\lambda^{r-1}}
\,\hbox{\vrule height 14pt depth 8pt width 0.5pt}_{\,\lambda=0}.
\mytag{3.11}
$$
The equality $W\cdot U=1$ derived from \mythetag{3.9} means that the 
matrix \mythetag{3.10} with the components \mythetag{3.11} is inverse
to the Vandermonde matrix \mythetag{3.5}.\par
     Since $(W^{\sssize\top})^{-1}=U^{\sssize\top}$, we can apply the transpose 
of the matrix $U$ in order to solve the matrix equation \mythetag{3.4}.
Its solution is written as
$$
\hskip -2em
a_q=U_{nq}=\frac{1}{(n-1)!}\,\frac{d^{\kern 1pt n-1}P_q(\lambda)}{d\lambda^{n-1}}
\,\hbox{\vrule height 14pt depth 8pt width 0.5pt}_{\,\lambda=0}=
\frac{1}{\dsize\prod^n_{s\neq q}(\lambda_q-\lambda_s)}.
\mytag{3.12}
$$
Along with solving the matrix equation \mythetag{3.4}, the quantities 
\mythetag{3.12} solve the polynomial equation \mythetag{3.3} as well. 
\par
     Now, using the quantities \mythetag{3.12} as coefficients, we define the 
following linear combination of the exponential functions \mythetag{1.1}:
$$
\hskip -2em
f_n(x)=a_1\,e^{\kern 0.5pt\lambda_1x}+\ldots+a_n\,e^{\kern 0.5pt\lambda_nx}.
\mytag{3.13}
$$
Taking into account \mythetag{3.1}, \mythetag{3.2} and \mythetag{3.3}, 
for $f_n(x)$ we derive the power expansion
$$
\hskip -2em
f_n(x)=\frac{x^{n-1}}{(n-1)!}+\sum^\infty_{q=1}\frac{B_{n\,q}\,x^{n-1+q}}{(n-1+q)!}.
\mytag{3.14}
$$
The coefficients $B_{n\,q}$ in \mythetag{3.14} are given by the formula
$$
\hskip -2em
B_{n\,q}=\sum^n_{i=1}a_i\,\lambda_{\kern 1pt i}^{n-1+q}.
\mytag{3.15}
$$
\mylemma{3.1} For mutually distinct numbers $\lambda_{\kern 1pt i}\neq\lambda_j$ and 
for $q\geqslant 1$ the quantities $B_{n\,q}$ from \mythetag{3.15} obey the
recurrent relationship  
$$
\hskip -2em
B_{n+1\,q}=B_{n\,q}+\lambda_{n+1}\,B_{n+1\,q-1}.
\mytag{3.16}
$$
\endproclaim
     Lemma~\mythelemma{3.1} is easily proved by means of direct 
calculations using \mythetag{3.12} and \mythetag{3.15}. 
Substituting $q=0$ \pagebreak into \mythetag{3.15} and taking into 
account the matrix equality \mythetag{3.4} for the quantities 
$a_1,\,\ldots,\,a_n$, we find that 
$$
\hskip -2em
B_{n\,0}=1\text{\ \ for all \ }n\geqslant 1.
\mytag{3.17}
$$ 
If $n=1$, the formula \mythetag{3.12} turns to $a_q=1$. Then \mythetag{3.15}
yields
$$
\hskip -2em
B_{1\,q}=\lambda_1^q\text{\ \ for all \ }q\geqslant 0.
\mytag{3.18}
$$ 
The formula \mythetag{3.18} can be derived directly from the first Tailor 
expansion \mythetag{3.1}.\par
     The formulas \mythetag{3.17} and \mythetag{3.18} along with the 
recurrent relationship \mythetag{3.16} are sufficient to determine 
all of $B_{n\,q}$ inductively on two parameters $n$ and $q$. 
\mylemma{3.2} The quantities $B_{n\,q}$ defined for mutually distinct 
lambdas $\lambda_{\kern 1pt i}\neq\lambda_j$ through the formulas \mythetag{3.12} 
and \mythetag{3.15} are given by the explicit formula
$$
\hskip -2em
B_{n\,q}\ =\kern -0.7em\msum{{\phantom{n}}}_{1\leqslant i_1\leqslant
\,\ldots\,\leqslant i_q\leqslant n}
\kern -1em\lambda_{i_1}\kern -1pt\cdot\ldots\cdot\lambda_{i_q}
\text{, \ where \ }n\geqslant 1\text{\ \ and \ }
q\geqslant 1. 
\mytag{3.19}
$$
\endproclaim
     In order to prove Lemma~\mythelemma{3.2} it is sufficient to verify
the formulas \mythetag{3.16} and \mythetag{3.18} upon substituting 
\mythetag{3.19} into them. The formula \mythetag{3.19}, along with the
formula \mythetag{3.17}, determines all of the quantities $B_{n\,q}$. 
\par
     Let's denote through $N_{n\,q}$ the number of summands in the formula
\mythetag{3.19}. This number is estimated in the following way:
$$
\hskip -2em
N_{n\,q}\leqslant n^q.
\mytag{3.20}
$$
Let's recall that in present section we consider the special case where the number
of clusters $m=1$, $k_1=n$, and $\Lambda_1=0$. Therefore let's denote 
$$
\hskip -2em
\varepsilon=\max(|\lambda_1|,\ldots,|\lambda_n|).
\mytag{3.21}
$$
The notation \mythetag{3.21} is a version of \mythetag{1.6} adapted to our present 
case. Applying \mythetag{3.20} and \mythetag{3.21} to the summands in 
\mythetag{3.14}, for $q\geqslant 1$ we get
$$
\hskip -2em
\left|\frac{B_{n\,q}\,x^{n-1+q}}{(n-1+q)!}\right|\leqslant
\frac{n^q\,\varepsilon^q\,|x|^{n-1+q}}{(n-1)!\,n\,(n+1)\cdot\ldots\cdot(h+q-1)}.
\mytag{3.22}
$$
A weaker estimate is sufficient for our purposes. Therefore from \mythetag{3.22}
we derive
$$
\hskip -2em
\left|\frac{B_{n\,q}\,x^{n-1+q}}{(n-1+q)!}\right|\leqslant
\frac{\varepsilon^q\,|x|^{n-1+q}}{(n-1)!}\text{, \ where \ }q\geqslant 1.
\mytag{3.23}
$$
From \mythetag{3.23} one can easily derive the norm estimate
$$
\hskip -2em
\left\Vert\frac{B_{n\,q}\,x^{n-1+q}}{(n-1+q)!}\right\Vert\leqslant
\sqrt{\frac{2\,\pi}{2\,n+2\,q-1}}\,\frac{\pi^{n-1}}{(n-1)!}\,
(\pi\,\varepsilon)^q.
\mytag{3.24}
$$
in term of the $L^2$-norm of the Hilbert space $\Cal H=L^2([-\pi,+\pi])$. 
Since $n=\const$ in \mythetag{3.1} and \mythetag{3.24}, the estimate
\mythetag{3.24} can be simplified. \pagebreak For this purpose we introduce 
the following constant that does not depend on $q\geqslant 1$:
$$
\hskip -2em
C_n=\sqrt{\frac{2\,\pi}{2\,n+1}}\,\frac{\pi^{n-1}}{(n-1)!}.
\mytag{3.25}
$$
Using the constant \mythetag{3.25}, the estimate \mythetag{3.24} is 
simplified as
$$
\hskip -2em
\left\Vert\frac{B_{n\,q}\,x^{n-1+q}}{(n-1+q)!}\right\Vert\leqslant
C_n\,(\pi\,\varepsilon)^q\text{, \ where \ }q\geqslant 1.
\mytag{3.26}
$$
If $\pi\,\varepsilon<1/2$, then the estimate \mythetag{3.26} produces
an estimate for the function $f_n(x)$ from \mythetag{3.13} and
\mythetag{3.14}. Here is this estimate:
$$
\hskip -2em
\left\Vert f_n(x)-\frac{x^{n-1}}{(n-1)!}\right\Vert\leqslant
\frac{C_n\,\pi\,\varepsilon}{1-\pi\,\varepsilon}\leqslant
2\,C_n\,\pi\,\varepsilon.
\mytag{3.27}
$$
\par
     Due to \mythetag{3.13} the function $f_n(x)$ is a linear combination 
of the exponential functions $e^{\kern 0.5pt\lambda_1x},\,\ldots,
\,e^{\kern 0.5pt\lambda_nx}$, i\.\,e\. 
$f_n(x)\in\Span(e^{\kern 0.5pt\lambda_1x},\,\ldots,e^{\kern 0.5pt\lambda_nx})$.
Therefore we can formulate the following theorem.
\mytheorem{3.1} For any $n$ mutually distinct complex quantities $\lambda_1,\,
\ldots,\,\lambda_n$ tending to zero there is a function $f(x)$ belonging
to the subspace 
$$
L=\Span(e^{\kern 0.5pt\lambda_1x},\,\ldots,e^{\kern 0.5pt\lambda_nx}).
$$ 
of the Hilbert space of square integrable functions $\Cal H=L^2([-\pi,+\pi])$
and such that
$$
\hskip -2em
\Vert f(x)-x^{n-1}\Vert\to 0
\mytag{3.28}
$$
as $\lambda_1,\,\ldots,\,\lambda_n$ tend to zero.
\endproclaim
     Theorem~\mythetheorem{3.1} is immediate from the inequality 
\mythetag{3.27}. It is important to note that the norm convergence 
in \mythetag{3.28} is irrespective of any mutual relations of 
$\lambda_1,\,\ldots,\,\lambda_n$ and is irrespective of the individual 
convergence rates of $\lambda_{\kern 1pt i}\to 0$. 
\par
     Having $n$ mutually distinct complex quantities $\lambda_1,\,\ldots,
\,\lambda_n$ converging to zero, we can take a part of them $\lambda_1,
\,\ldots,\,\lambda_s$, where $1\leqslant s\leqslant n$. Applying the 
theorem~\mythetheorem{3.1} to all of such parts, we derive the following result. 
\mytheorem{3.2} For any $n$ mutually distinct complex quantities $\lambda_1,\,
\ldots,\,\lambda_n$ tending to zero there are $n$ function $f_1(x),\,
\ldots,\,f_n(x)$ belonging to the subspace 
$$
L=\Span(e^{\kern 0.5pt\lambda_1x},\,\ldots,e^{\kern 0.5pt\lambda_nx}).
$$ 
of the Hilbert space of square integrable functions $\Cal H=L^2([-\pi,+\pi])$
and such that
$$
\hskip -2em
\Vert f_s(x)-x^{s-1}\Vert\to 0
\mytag{3.29}
$$
as $\lambda_1,\,\ldots,\,\lambda_n$ tend to zero for all $s=1,\,\ldots,\,n$.
\endproclaim
     Again, it is important to note that the norm convergences \pagebreak
in \mythetag{3.29} are irrespective of any mutual relations of $\lambda_1,\,
\ldots,\,\lambda_n$ and are irrespective of the individual convergence rates 
of $\lambda_{\kern 1pt i}\to 0$.\par
     The quantities $\lambda_1,\,\ldots,\,\lambda_n$ tending to zero can go
through some discrete sets of values or their convergence can be a continuous
process. In both cases, relying on the definitions~\mythedefinition{2.1} 
and \mythedefinition{2.2}, we can derive the following theorem.
\mytheorem{3.3} For any $n$ mutually distinct complex quantities $\lambda_1,\,
\ldots,\,\lambda_n$ tending to zero the span of exponential functions
$$
L=\Span(e^{\kern 0.5pt\lambda_1x},\,\ldots,e^{\kern 0.5pt\lambda_nx})
$$ 
converges to the span of polynomials 
$$
M=\Span(1,\,x,\,\ldots,\,x^{n-1})
$$
in the Hilbert space of square integrable functions $\Cal H=L^2([-\pi,+\pi])$. 
\endproclaim
     Theorem~\mythetheorem{3.3} is immediate from the previous 
theorem~\mythetheorem{3.2}.
\head
4. The case of multiple clusters. 
\endhead     
     Now we proceed to the general case where lambdas are subdivided into
$m$ clusters with $k_1,\,\ldots,\,k_m$ being the multiplicities of clusters
(see \mythetag{1.3}). They tend to $m$ mutually distinct complex numbers
$\Lambda_1,\,\ldots,\,\Lambda_m$ according to \mythetag{1.5} or \mythetag{1.7}. 
Therefore in this case we introduce the deflection numbers
$$
\delta_{ij}=\lambda_{ij}-\Lambda_i
$$
tending to zero and, instead of \mythetag{3.1}, we write 
$$
\aligned
&e^{\kern 0.5pt\lambda_{i1}x}=e^{\Lambda_ix}\left(\!1+\delta_{i1}\,x+\ldots
+\frac{\delta_{i1}^{n-1}\,x^{n-1}}{(n-1)!}\right)+\ldots
=\sum^\infty_{q=1}\frac{\delta_{i1}^q\,x^q\,e^{\Lambda_ix}}{q!},\\
\vspace{-1ex}
&.\ .\ .\ .\ .\ .\ .\ .\ .\ .\ .\ .\ .\ .\ .\ .\ .\ .\ .\ .\ .\ .\ .\ .\ 
.\ .\ .\ .\ .\ .\ .\ .\ .\ .\ .\ .\ .\ .\ 
\\
\vspace{-0.4ex}
&e^{\kern 0.5pt\lambda_{ik_i}x}=e^{\Lambda_ix}\left(1+\delta_{ik_i}\,x+\ldots
+\frac{\delta_{ik_i}^{n-1}x^{n-1}}{(n-1)!}\right)+\ldots
=\sum^\infty_{q=1}\frac{\delta_{ik_i}^q\,x^q\,e^{\Lambda_ix}}{q!}. 
\endaligned\quad
$$
As a result, instead of Theorem~\mythetheorem{3.1}, here we get the 
following theorem. 
\mytheorem{4.1} For a set of mutually distinct complex quantities 
$\lambda_{\kern 0.5pt ij}$, where\linebreak $i=1,\,\ldots,m$ and $j=1,\,\ldots,k_i$,
tending to $m$ mutually distinct complex numbers $\Lambda_1,\,\ldots,\,
\Lambda_m$ so that $\lambda_{ij}\to\Lambda_i$ there is a function $f_s(x)$ 
belonging to the subspace 
$$
L=\Span(\{e^{\kern 0.5pt\lambda_{ij}x}\}^{i=1,\,\ldots,\,m}_{j=1,\,\ldots,\,k_i})
$$ 
of the Hilbert space of square integrable functions $\Cal H=L^2([-\pi,+\pi])$
and such that
$$
\hskip -2em
\Vert f_s(x)-x^{k_s-1}\,e^{\Lambda_sx}\Vert\to 0\text{\ \ as \ }
\lambda_{ij}\to\Lambda_i.
\mytag{4.1}
$$
\endproclaim
     The proof of Theorem~\mythetheorem{4.1} is basically the same as the proof
of Theorem~\mythetheorem{3.1}. The main difference is that, instead of
the polynomial $x^{n-1}$, \pagebreak here in Theorem~\mythetheorem{4.1} we have 
the expo-polynomial $x^{k_i-1}\,e^{\Lambda_ix}$.\par
     The expo-polynomial $x^{k_i-1}\,e^{\Lambda_ix}$ is not unique. Repeating
the reasons used in deriving Theorem~\mythetheorem{3.2} from 
Theorem~\mythetheorem{3.1}, we can write the following theorem that provides
as many expo-polynomials as the exponential functions we have. 
\mytheorem{4.2} For a set of mutually distinct complex quantities 
$\lambda_{\kern 0.5pt ij}$, where\linebreak $i=1,\,\ldots,m$ and where 
$j=1,\,\ldots,k_i$, tending to $m$ mutually distinct complex numbers 
$\Lambda_1,\,\ldots,\,\Lambda_m$ so that $\lambda_{ij}\to\Lambda_i$ there 
is a set of functions $f_{ij}(x)$, where\linebreak $i=1,\,\ldots,m$ and 
$j=1,\,\ldots,k_i$, belonging to 
the subspace 
$$
L=\Span(\{e^{\kern 0.5pt\lambda_{ij}x}\}^{i=1,\,\ldots,\,m}_{j=1,\,\ldots,\,k_i})
$$ 
of the Hilbert space of square integrable functions $\Cal H=L^2([-\pi,+\pi])$ and 
such that
$$
\hskip -2em
\Vert f_{ij}(x)-x^{\kern 1pt j-1}\,e^{\Lambda_ix}\Vert\to 0\text{\ \ as \ }
\lambda_{ij}\to\Lambda_i.
\mytag{4.2}
$$
\endproclaim
     It is important to note that the norm convergences in \mythetag{4.1} and
\mythetag{4.2} are irrespective of any mutual relations between $\lambda_{ij}$ 
and are irrespective of the individual convergence rates of 
$\lambda_{ij}\to\Lambda_i$.\par
     The quantities $\lambda_{ij}$ tending to $\Lambda_i$ within their clusters
can go through some discrete sets of values or their convergence can be a continuous
process. In both cases, relying on the definitions~\mythedefinition{2.1} 
and \mythedefinition{2.2}, we can derive the following theorem.
\mytheorem{4.3} For a set of mutually distinct complex quantities 
$\lambda_{\kern 0.5pt ij}$, where\linebreak $i=1,\,\ldots,m$ and $j=1,\,\ldots,k_i$,
tending to $m$ mutually distinct complex numbers $\Lambda_1,\,\ldots,\,
\Lambda_m$ so that $\lambda_{ij}\to\Lambda_i$ the span of exponential functions
$$
L=\Span(\{e^{\kern 0.5pt\lambda_{ij}x}\}^{i=1,\,\ldots,\,m}_{j=1,\,\ldots,\,k_i})
$$ 
converges to the span of expo-polynomials 
$$
M=\Span(\{x^{\kern 1pt j-1}\,e^{\Lambda_ix}\}^{i=1,\,\ldots,\,m}_{j=1,\,\ldots,\,k_i})
$$
in the Hilbert space of square integrable functions $\Cal H=L^2([-\pi,+\pi])$. 
\endproclaim
     Theorem~\mythetheorem{4.3} is the main result of this paper. It is immediate
from Theorem~\mythetheorem{4.2}.
\head
5. Acknowledgment.
\endhead
     I am grateful to my friend A.~S.~Vishnevskiy who stimulated my interest to 
a mathematical problem originated in his business. This paper is the first step in 
my attempts toward solving this problem.\par

\enddocument
\end